\documentclass[12pt]{amsart}

\begin{document}
\title[MOORE-PENROSE INVERSE AND EP ELEMENTS ]{ ON THE MOORE-PENROSE INVERSE,\\
EP BANACH SPACE OPERATORS,\\
AND EP BANACH ALGEBRA ELEMENTS }
\author[ENRICO BOASSO ]{ENRICO BOASSO}

\begin{abstract}
The main concern of this note is the Moore-Penrose inverse in the
context of Banach spaces and algebras. Especially attention will be given
to a particular class of elements with the aforementioned inverse, namely EP  
Banach space operators and Banach algebra elements, 
which will be studied and characterized extending well-known results obtained
in the frame of Hilbert space operators and $C^*$-algebra elements.\par
\par
\end{abstract}
\maketitle

\noindent \bf{1. Introduction}\rm \vskip.3cm

\indent The Moore-Penrose inverse is a notion that was introduced for 
matrices, see [21], and whose development has produced a wide literature. 
In the presence of an involution, in a Hilbert space or in a $C^*$-algebra,
it is natural to extend and study the main properties 
of the aforesaid notion, see for example [4], [6], [8], [11], [12], 
[15], [18], and [19].\par

\indent With the extension to general Banach space operators and more 
generally Banach algebra elements of the concept of hermitian element,
see [2], [3], [5], [9], [20], [25], and [26], came also an extended notion of
Moore-Penrose inverse, due to V. Rakocevic, see [22] and [23].\par

\indent In particular, the class of $C^*$-algebra elements which commute with their
Moore-Penrose inverse, the so-called EP elements, has been of
especial interest. In fact, several results characterizing these
elements and when the product of two EP elements
is again EP have been obtained, see [1], [6], [8], [13], [14], and [15].\par

\indent In this work the Moore-Penrose inverse in Banach spaces and 
algebras will be studied. In first place, several basic facts, as the
relationships between the Moore-Penrose inverse and, on the one hand, closed
invariant subspaces and, on the other, the adjoint of an operator will be
considered. In addition, several characterizations of the
aforementioned inverse will be proved.\par

\indent On the other hand, EP Banach space operators and
Banach algebra elements will be studied. In fact, 
they will be characterized extending to this context
well-known characterizations of Hilbert space operators
and of $C^*$-algebra elements, see [6], [12], and [15].
Furthermore, the problem of characterizing when
the product of two EP elements is again EP will
be considered extending results of [13], [8],
[14], and [15], see also [1].\par

\indent The work is organized as follows. In section 2
several preliminary definitions and results will be
recalled. In section 3 the Moore-Penrose inverse
in Banach spaces and algebras will be studied.
In section 4 EP elements will be characterized,
and in section 5 it will be considered the problem of
determining when the product of two EP elements is again EP.\par

\vskip.3cm 
\noindent \bf{2. Preliminary Definitions and Results}\rm \vskip.3cm
\indent From now on, $X$ will denote a Banach space, and $L(X, Y)$ the
Banach algebra of all bounded and linear maps defined on $X$ with values 
in the Banach space $Y$. As usual, when $X=Y$, $L(X,Y)$ will be
denoted by $L(X)$, the Banach algebra of all operators with domain $X$.
In addition, if $T\in L(X)$, then $N(T)$ and $R(T)$ will stand for the null space and the
range of $T$, respectively. \par
 
\indent On the other hand, $A$ will denote a unital Banach algebra, that is
a Banach algebra with a unit element $e$ such that $\parallel e\parallel =1$. 
If $a\in A$, then $L_a \colon A\to A$
and $R_a\colon A\to A$ will denote the map defined by
left and right multiplication, respectively:
$$
L_a(x)=ax, \hskip2truecm R_a(x)=xa,
$$   
where $x\in A$. Moreover, the following notation will be used:

\begin{align*}
&N(L_a)= a^{-1}(0),& &R(L_a)=aA,\\
&N(R_a)= a_{-1}(0),& &R(R_a)=Aa.\\
\end{align*}

\indent Recall that an element $a\in A$ is called \it{regular}, \rm 
if it has a \it{generalized inverse}, \rm namely if there exists $b\in A$ such that
$$
a=aba.
$$

\indent Furthermore, a generalized inverse $b$ of a regular
element $a\in A$ will be called \it{normalized}, \rm if $b$ is regular
and $a$ is a generalized inverse of $b$, equivalently,
$$
a=aba, \hskip2truecm b=bab,
$$
see for example [11], [12], and [18].  \par
\indent Note that if $b$ is a generalized inverse of $a$,
then $c=bab$ is a normalized generalized inverse
of $a$.\par

\indent Next follows the key notion in the definition of the 
Moore-Penrose inverse in context of Banach spaces and algebras.\par

\newtheorem*{def1}{Definition 1} 
\begin{def1}Given a unital Banach algebra $A$, an element $a\in A$ is said to be hermitian,
if $\parallel exp(ita)\parallel =1$,  for all $ t\in\Bbb R$.
\end{def1}

\indent As regard equivalent definitions and the main properties of hermitian Banach 
algebra elements and Banach space operators, see for example [2], [3], [20], 
[9], [26], [25], and [5]. In the following remark some relevant facts
will be considered.\par

\newtheorem*{rem2}{Remark 2} 
\begin{rem2}\rm In the conditions of the previous 
definition, recall that if $A$ is a $C^*$-algebra, then
$a\in A$ is hermitian if and only if $a$ is self-adjoint,
see Proposition 20, section 12, Chapter I of [5].
Furthermore, $\mathcal{H}=\{a\in A\colon \hbox{  }$a$ \hbox{ is 
hermitian}\}\subseteq A$ is a closed linear vector space over the real field,
see [26] and Chapter 4 of [9]. Since $A$ is unital, $e\in \mathcal{H}$,
which implies that $a\in \mathcal{H}$ if and only if $e-a\in \mathcal{H}$.\par
\indent Observe also that necessary and sufficient for $a\in A$ to be hermitian
is that $\parallel exp(ita)\parallel\le 1$,
$\forall$ $t\in \Bbb R$. In fact, if this condition holds, then
$$
1=\parallel exp(ita)exp(-ita)\parallel\le \parallel exp(ita)\parallel\parallel exp(-ita)\parallel\le 1.
$$
\indent On the other hand, recall that there is an equivalent definition of hermitian element in terms of numerical range,
see [17], [20], Chapter 4 of [9], section 10, Chapter I of [5], and also [2].\par 
\indent When $A=L(X)$, $X$ a Banach space, the hermitian idempotents
will be of particular interest. Following [3], the set of all these elements will be
denoted by $\mathcal{E}(X)$, that is 
$$
 \mathcal{E}(X)=\{ P\in L(X)\colon P^2=P, \hbox{ and }P \hbox{ is a hermitian operator}\}.
$$
\noindent Recall that, according to Theorem 2.2 of [20], $P=Q$ if and only 
if $R(P)=R(Q)$, where $P$ and $Q$ belong to
$ \mathcal{E}(X)$.\par
\indent In addition, following the Definitions of page 108 and 110 in [3], $\mathcal{M}(X)$
will denote the collection of all subspaces $M$ of $X$ such that there is a (necessarily
unique) $P\in \mathcal{E}(X)$ with $M=R(P)$. In this case, for $M\in \mathcal{M}(X)$,
the unique $P\in \mathcal{E}(X)$ 
such that $R(P)=M$ will be denoted by $P_M$. Further, $N(P_M)$ will be denoted by
$M'$.\par
\end{rem2}

\indent Now the central notion of the present note will be considered.\par

\indent Let $A$ be a Banach algebra and $a\in A$. The main concern of
this work consists in the elements $a$ for which there exists $x\in A$
satisfying the following conditions:

\begin{align*}
 & \rm{ (i)}\hskip.23cm  \it a= axa,\hskip2truecm x=xax,\\
 & \rm { (ii)}\hskip.15cm \it ax \rm \hbox{ and } \it xa \rm \hbox{ are hermitian}.
\end{align*}

\indent In Lemma 2.1 of [22] it was proved that given $a\in A$, there exists at most
one $x$ such that the previous conditions hold. This fact led to the following
definition\par
 
\newtheorem*{def3}{Definition 3} 
\begin{def3} Let $A$ be a unital Banach algebra and $a\in A$. If there exists 
$x\in A$ such that the previous conditions are satisfied, the element $x$ 
will be called the Moore-Penrose inverse of $a$, and it will be
denoted by $a^{\dag}$.
\end{def3}

\indent As regard the Moore-Penrose inverse in Banach algebra, 
see [22], where this concept was introduced, [23],
the continuation of [22], and see also [19],
where Moore-Penrose inverse of Banach space operators were
considered. For the original definition of the Moore-Penrose
inverse for matrices, see [21].\par
 
\newtheorem*{rem4}{Remark 4} 
\begin{rem4}\rm In the same conditions of Definition 3, it is clear that if $a\in A$ has a Moore-Penrose inverse,
then $a^{\dag}$ also has one and $(a^{\dag})^{\dag}=a$.
Furthermore, invertible elements and hermitian idempotents have a Moore-Penrose inverse.\par

\indent Note that the norm is a fundamental notion involved in Definitons 1 and 3. Actually,
even when the underlying space is the same and two equivalent
norms are considered, it is possible that an operator is hermitian (resp. has a Moore-Penrose
inverse) in one of the norms but not in the other.\par
\indent In fact, define $T(x,y) =(\frac{1}{2}x -\frac{1}{2}y, -\frac{1}{2}x+\frac{1}{2}y)$, where $(x,y)\in \Bbb C^2$. Then, 
it is not difficult to prove that, on the one hand, $T$ is a hermitian idempotent with 
the Euclidean norm, while on the other, $T$ is not hermitian in 
$(\Bbb C, \parallel \cdot\parallel_1)$, where $\parallel (x,y)\parallel_1 = |\hbox{ }x\hbox{ }| +|\hbox{ } y\hbox{ }|$, $(x,y)\in C^2$.\par
\indent In addition, if $S(x,y)=(x-y, 0)$, $(x,y)\in \Bbb C^2$, then, while in the Hilbert space norm
$S^{\dag}$ exists, a rather painful
calculation shows that $S$ does not have a Moore-Penrose inverse in 
$(\Bbb C, \parallel \cdot\parallel_1)$.\par  

\indent On the other hand, recall that according to
Theorem 6 of [13], in a $\Bbb C^*$-algebra
the set of the elements with a Moore-Penrose inverse
coincide with the one of the regular elements.
However, in a Banach algebra this result does not hold in general
any more. In fact, $S\in A=L((\Bbb C, \parallel \cdot\parallel_1) )$
is an idempotent without a Moore-Penrose inverse. Moreover, in Theorem
4.6 of [19] it was proved that necessary and sufficient for
a Banach space of dimension greater than 3 to be a Hilbert space is that the set of all 
regular operators coincides with the one of all bounded and linear maps with a
Moore-Penrose inverse.\par

\end{rem4}

\indent In the following section the Moore-Penrose inverse in
Banach spaces and algebras will be studied.\par

\noindent \bf{3. Basic Properties of the Moore-Penrose Inverse }\rm \vskip.3cm

\indent In this section several characterizations concerning the Moore-Penrose inverse
in Banach spaces and algebras will be given. In particular, the relationships
between the aforesaid concept and, on the one hand,
 the adjoint of a Banach space operator and, on the other, closed invariant
subspaces will be considered.\par

\indent  For $T\in L(X)$,
the definitions of the Moore-Penrose inverse as an operator and as
an element of the algebra $A=L(X)$ coincide. On the other hand, 
the starting point of 
the next theorem will be the Banach algebra setting, in which
the notion under consideration will be related to
the Moore-Penrose inverse of Banach space operators.\par

\newtheorem*{theo5}{Theorem 5} 
\begin{theo5} Let $A$ be a unital Banach algebra, and $a\in A$.\par
(i) An element $b\in A$ is a normalized generalized inverse of $a$ if and only if
$L_b$ (resp. $R_b$) is a normalized generalized inverse of $L_a$ (resp. $R_a$) in $L(A)$.\par
(ii) For a regular element $a\in A$,
necessary and suffcient for $b\in A$ to be the Moore-Penrose inverse of $a$ in $A$
is that $L_b$ (resp. $R_b$) is the Moore-Penrose inverse of $L_a$ (resp. $R_a$) in $L(A)$. Moreover, in this case
$(L_a)^{\dag}=L_{a^{\dag}}$ (resp. $(R_a)^{\dag}=R_{a^{\dag}}$).\par
\end{theo5}  
\begin{proof} The first statement is clear.\par
\indent As regard the second statement, a straightforward calculation
shows that
$$
exp(itL_d)=L_{exp(itd)},
$$
where $d\in A$, and $t\in \Bbb{R}$. However, since $\parallel e\parallel =1$,
$\parallel L_d\parallel =\parallel d\parallel$. Therefore, for $t\in \Bbb{R}$,
$$
\parallel exp(itL_d)\parallel =\parallel L_{exp(itd)}\parallel =\parallel exp(itd)\parallel.
$$
\indent Next consider $b\in A$, a normalized generalized inverse of $a$.
Since $L_{ab} =L_aL_b$ and $L_{ba} =L_bL_a$, the above equality
shows that $b$ is the Moore-Penrose inverse of $a$ if and only if $L_b$
is the Moore-Penrose inverse of $L_a$.\par
\indent In a similar way, it can be proved the second statement considering
$R_a$ instead of $L_a$.\par 
\end{proof}

\indent It is well known that necessary and sufficient for $T\in L(X)$ to be 
regular is that $N(T)$ and $R(T)$ are closed and complemented subspaces
of $X$, see Theorem 1 of [7] or Theorem 3.8.2 of [10]. In the following
theorem, the corresponding characterization for Banach space 
operators and Banach algebra elements with a Moore-Penrose inverse
will be given.\par

\newtheorem*{theo6}{Theorem 6} 
\begin{theo6} Let $X$ be a Banach space and $T\in L(X)$.
Then, the following statements are equivalent:\par
\noindent  (i) \hskip.18cm $T$ has a Moore-Penrose inverse,\par
\noindent (ii) there exist P and Q in $\mathcal{ E}(X)$,
such that $N(P)=N(T)$ and $R(Q)=R(T)$.\par
\indent Let $A$ be a unital Banach algebra, and consider an element $a\in A$. Then, the following
statements are equivalent:\par
\noindent  (i) \hskip.18cm The map $L_a$ (resp. $R_a)$ has a Moore-Penrose inverse in $L(A)$,\par
\noindent (ii) there exist P and Q in $\mathcal{ E}(A)$, such that
$N(P)=a^{-1}(0)$ and $R(Q)=aA$ (resp. $N(P)=a_{-1}(0)$ and $R(Q)=Aa$.\par
\indent Furthermore, if such $P$ and $Q$ exist, then they are unique.\par 
 \end{theo6}
\begin{proof} 
\indent If $T\in L(X)$ has a Moore-Penrose inverse $T^{\dag}$, then consider
$P= T^{\dag}T$ and $Q=TT^{\dag}$.\par
\indent In order to verify the converse implication, consider the
invertible operator $T'\in L(R(P), R(T))$, 
$$
T'=T\mid_{R(P)}^{R(T)}\colon R(P)\to R(T),
$$ 
and define $S\in L(X)$ as follows:
$$
S\mid N(Q)\equiv 0,\hskip1cm S\mid_{R(T)}^{R(P)} =(T')^{-1}\in L(R(T),R(P)).
$$
\indent A straightforward calculation proves that $S$ is a normalized
generalized inverse of $T$. On the other hand, since $TS$ and $Q$
are idempotents of $L(X)$ such that
$$
R(TS)=R(T)=R(Q), \hskip1cm N(TS)=N(S)=N(Q),
$$
it is clear that $TS=Q$. In particular, $TS\in \mathcal{E}(X)$. Similarly,
$ST\in \mathcal{E}(X)$. Therefore, $S$ is the Moore-Penrose inverse of $T$.\par
\indent Moreover, if $P'$ and $Q'$ are two other hermitian idempotents
that satisfy the above conditions, then $R(Q) =R(Q')$ and $R(I-P)=R(I-P')$.
Then, according to Theorem 2.2 of [20], or to Remark 2, $P=P'$
and $Q=Q'$.\par

\indent The rest of the Theorem can be deduced from what has been proved.\par
 
\end{proof} 
\indent Next follows the relationships between the Moore-Penrose inverse
of a Banach space operator and a closed invariant subspace.\par

\newtheorem*{prop7}{Proposition 7} 
\begin{prop7}Let $X$ be a Banach space, and 
$T\in L(X)$ such that there exists $T^{\dag}$. Let $Y$ be a closed
subspace of $X$ invariant  under $T$ and $T^{\dag}$.\par
(i) \hskip.2cm If $T'$ (resp. $T^{\dag'}$) denotes the operator induced by $T$ (resp. by $T^{\dag}$)
on $Y$, then $T^{\dag'}$ is the Moore-Penrose inverse of $T'$ on $Y$.\par
(ii) If $T_Y$ (resp. $T^{\dag}_Y$) denotes the operator induced by $T$
(resp. by $T^{\dag}$) on the quotient Banach space $X/Y$, then
$T^{\dag}_Y$ is the Moore-Penrose inverse of $T_Y$ on $X/Y$.\par
\indent On the other hand, consider $X_i$, $i=1$, $2$, two Banach spaces, and $T_i\in L(X_i)$, $i=1$, $2$,
two operators such that there exist $T^{\dag}_i\in L(X_i)$, $i=1$, $2$. Define 
the Banach space $X=X_1\oplus X_2$, with the norm $\parallel x_1\oplus x_2\parallel=
max\{ \parallel x_1\parallel, \parallel x_2\parallel\}$. Then $T^{\dag}_1\oplus T^{\dag}_2$ 
is the Moore-Penrose inverse of $T=T_1\oplus T_2$ on $X$.\par
\end{prop7}  
\begin{proof} it is clear that $T^{\dag'}$ is a normalized
generalized inverse of $T' $. Moreover, according to the first
statement of Proposition 4.12 of [9], $T^{\dag'}$ is the 
Moore-Penrose inverse of $T' $.\par

\indent The second statement can be proved in a similar way,
using the statement (ii) of Proposition 4.12 of [9].\par

\indent As regard the last part of the Proposition, an easy
calculation shows that $S=T^{\dag}_1\oplus T^{\dag}_2$ is a normalized
generalized inverse of $T$. In addition, given Banach space operators
$S_i\in L(X_i)$, $i=1$, $2$, it is not difficult
to prove that,  
$$
exp(S_1\oplus S_2)=exp(S_1)\oplus exp(S_2),
$$
which implies that $\parallel exp(S_1\oplus S_2)\parallel\le
max\{\parallel exp(S_1)\parallel, \parallel exp(S_2)\parallel\}$.
In particular, if $S_i=itT_iT^{\dag}_i$, $i=1$, $2$,
then $\parallel exp(itTS)\parallel\le 1$,
where $t\in \Bbb R$.\par

\indent Similarly, $\parallel exp(itST)\parallel\le 1$.
Now well, according to Remark 2, $S= T^{\dag}_1\oplus T^{\dag}_2$ is
the Moore-Penrose inverse of $T$.\par 
\end{proof}

\indent In the setting of Banach algebras,
Proposition 7 can be generlized in the following way.\par

\newtheorem*{prop8}{Proposition 8} 
\begin{prop8} Let $A$ be a unital Banach algebra, and $a\in A$ such
that there exists $a^{\dag}$.\par
(i)  \hskip.2cm If $B$ is a subalgebra of $A$ such that $a$ and $a^{\dag}$
belong to $B$, then $a^{\dag}$ is the Moore-Penrose inverse
of $a$ in $B$.\par
(ii) If $J$ is a proper and closed bilateral ideal of $A$, then $(\tilde{a})^{\dag}=\tilde{a^{\dag}}$, where if $d\in A$,
$\tilde{d}$ denotes the quotient class of $d$ in $A/J$.
\end{prop8} 
\begin{proof} The first statement is clear.\par

\indent An easy calculation shows that $\tilde{a^{\dag}}$ is a normalized
generalized inverse of $\tilde{a}$ in the unital Banach algebra $A/J$. \par

\indent Next consider the map $L_a\in L(A)$. It is clear that
$L_a(J)\subseteq J$. Moreover, following the notation of
Proposition 7, $L_{aJ}=L_ {\tilde{a}}\in L(A/J)$. 
Similarly, $L_{a^{\dag}J}=L_ {\tilde{a^{\dag}}}\in L(A/J)$. Now well,
according to the second statement of Theorem 5 and
Proposition 7, 
$$
(L_a)^{\dag}= L_{a^{\dag}}, \hskip1cm (L_{aJ})^{\dag}=L_{a^{\dag}J}.
$$
\indent Therefore,
$$
(L_ {\tilde{a}})^{\dag}=(L_{aJ})^{\dag}= L_{a^{\dag}J}=L_ {\tilde{a^{\dag}}}, 
$$
with which, according again to the second statement of Theorem 5,
the proof is concluded.\par 
\end{proof}

\indent In the last Theorem of the present section, the relationship between the Moore-Penrose
inverse and the adjoint in Banach spaces will be studied. First of all
some notation will be recalled. If $X$ is a Banach space, then its dual will be denoted
by $X^*$. In addition, if $T\in L(X)$, the adjoint map of $T$
will be denoted by $T^*$. Next follows a necessary remark.\par

\newtheorem*{rem9}{Remark 9} 
\begin{rem9}\rm Let $X$ and $Y$ be two Banach spaces and $U\in L(Y)$. Consider 
 $F\colon X\to Y$ an isometric isomorphism.
Then, an easy calculation shows that
$$
exp( F^{-1}UF)=F^{-1}exp(U)F.
$$
\indent In particular, $ F^{-1}UF$ is hermitian in $L(X)$ if and only if 
$U$ is hermitian in $L(Y)$.\par
\indent Furthermore, if $U$ has a Moore-Penrose inverse in $Y$,
then it is not difficult to prove that $F^{-1}U^{\dag}F$ is the Moore-Penrose
inverse of $F^{-1}UF$ in $X$.\par
 \end{rem9}

\newtheorem*{theo10}{Theorem 10} 
\begin{theo10} Let $X$ be a Banach space and $T\in L(X)$. \par
(i) \hskip.2cm If there exists $T^{\dag}$, then $T^*$
has a Moore-Penrose inverse and $(T^*)^{\dag}= (T^{\dag})^*$.\par
(ii) Suppose that there exist $(T^*)^{\dag}\in L(X^*)$ and $S\in L(X)$ such 
that $S^*= (T^*)^{\dag}$, then there exists $T^{\dag}$ and $S=T^{\dag}$.\par
\indent In particular, if $X$ is a reflexive space, then necessary and sufficient for $T\in L(X)$
to have a Moore-Penrose inverse is that $(T^*)^{\dag}$ exists.\par
\end{theo10}
\begin{proof} If $T^{\dag}$ exists, then $(T^{\dag})^*$
is a normalized generalized inverse of  $T^*$. Moreover, according
to Proposition 4.11 of [9], $(T^{\dag})^*$ is the Moore-Penrose
inverse of $T^*$.\par

\indent As regard the second part of the Proposition, according 
to what has been proved, there exists $((T^*)^*)^{\dag}$ and $((T^*)^*)^{\dag}=(S^*)^*\in
L((X^*)^*)$.\par

\indent Next consider $\tilde{T}$ (resp. $\tilde{S}$), the restriction of 
$(T^*)^*$ (resp. $(S^*)^*$) to the closed invariant subspace $\Lambda (X)\subseteq (X^*)^*$,
where
$$
\Lambda\colon X\to \Lambda (X)\subseteq (X^*)^*, \hskip.5cm \Lambda (x)(f)=f(x),\hskip.5cm (x\in X\hbox{ and }f\in X^*),
$$
is the canonical isometric identification of $X$
with a closed subspace of $(X^*)^*$. Then,
according to the first statement of Proposition 7, there exists 
$\tilde{T}^{\dag}$ and $\tilde{T}^{\dag}=\tilde{S}$. Now well,
according to Remark 9, there exists $T^{\dag}$ and 
$T^{\dag}=S$.

\indent The last statement is a consequence of the fact that, 
when $X$ is a reflexive space, $L(X^*)=\{S^*\colon S\in L(X)\}$.\par
\end{proof}

\noindent \bf{4. EP Banach Space Operators and Banach Algebra Elements  }\rm \vskip.3cm

\indent In this section a particular case of the elements with a Moore-Penrose inverse
will be considered, namely EP Banach space operators
and Banach algebra elements. As in the previous section, the basic properties of these objects
will be studied. Furthermore, several characterizations will be given extending well-known results 
obtained in the frame of Hilbert spaces and $C^*$-algebras, see 
[6], [12], and [15]. In first place, the definition 
of the aforementioned notion will be introduced.\par

\newtheorem*{def11}{Definition 11} 
\begin{def11}Given a unital Banach algebra $A$, $a\in A$ will be said an 
EP element, if there exists $a^{\dag}$, and  $aa^{\dag}= a^{\dag}a$. 
\end{def11}

\indent Observe that the name of the elements introduced in
Definition 11 derives fron the fact that 
the idempotents, projectors in the Banach space operator
context, $aa^{\dag}$ and $a^{\dag}a$ are equal. On the other hand,
In the following remark some of the basic results
regarding EP Banach algebra elements 
will be collected. Each of them can be 
proved with a direct argument.\par  

\newtheorem*{rem12}{Remark 12} 
\begin{rem12}\rm 
\indent In the same conditions of Definition 11, consider $a\in A$ such that $a^{\dag}$
exists. Observe that since $(a^{\dag})^{\dag}=a$, $a$
is EP if and only if $a^{\dag}$ is. 
In addition, a direct calculation proves that if $a\in A$
is EP, then so is $a^k$, for $k\ge 1$.
Moreover,
since according to Theorem 5, 
$(L_a)^{\dag}=L_{a^{\dag}}$, $a\in A$ is EP if and only if 
$L_a\in L(A)$ is. A similar equivalence can be obtained for the 
map $R_a\in L(A)$.\par
\indent Recall that the \it group inverse \rm of $a\in A$ is an element $b\in A$,
such that $a=aba$, $b=bab$, and $ab=ba$, see for example [24]. Note that
when $a$ is an EP element, then $a$ has a group inverse, which
coincides with the Moore-Penrose inverse.\par
\indent On the other hand, given a Banach space $X$ and $T\in L(X)$ and EP
operator, if $Y\subseteq X$ is a 
closed subspace of $X$ which is invariant
under $T$ and $T^{\dag}$, then
according to Theorem 7, $T'\in L(Y)$ and $T_Y\in L(X/Y)$ are
EP operators.\par 
\end{rem12}
 
\indent Next follows a characterization of EP Banach space operators
and Banach algebra elements.\par

\newtheorem*{prop13}{Proposition 13} 
\begin{prop13}Let $X$ be a Banach space and 
$T\in L(X)$. Then, the following statements are
equivalent:\par
\noindent (i) \hskip.2cm $T$ is an EP operator,\par
\noindent (ii) there exists $P\in \mathcal{E}(X)$ such that $N(P)=N(T)$
and $R(P)=R(T)$.\par
\indent Let $A$ be a unital Banach algebra, and consider a regular element $a\in A$ such that
$a^{\dag}$ exists.
Then, the following statements are equivalent:\par
\noindent (i)\hskip.22cm $a$ is EP,\par
\noindent (ii) there exists $P\in \mathcal{E}(A)$ such that 
$N(P)=a^{-1}(0)$ and $R(P)=aA,$ 
\noindent (iii) there exists $P\in \mathcal{E}(A)$ such that 
$N(P)=a_{-1}(0)$ and $R(P)=Aa.$ 
\indent Furthermore, if such $P$ exists, then it is unique.\par
\end{prop13}
\begin{proof} If $T$ is an EP operator, then the idempotent 
$P=TT^{\dag}= T^{\dag}T$ complies with the required property.\par

\indent On the other hand, according to Theorem 6, $T^{\dag}$ exists.
In addition, since $R(I-T^{\dag}T)=N(T)=N(P)=R(I-P)$, and since $I-P$
is an hermitian idempotent, see Remark 2, [26] or Chapter 4 of [9],
according to Theorem 2.2 of [20], $TT^{\dag}=P=T^{\dag}T$.\par

\indent Consider a Banach algebra $A$, and $a\in A$ a regular element such that
$a^{\dag}$ exists.
If $a$ is EP, then, according to Theorem 5 and the proof of
Theorem 6, $P=L_ {aa^{\dag}}$ satisfies the desired condition.\par
\indent In order to prove the converse implication, note that
according to Theorem 5, $L_a$ has a Moore-Penrose inverse in $L(A)$.
Moreover, thanks to what has been proved, $L_a\in L(A)$ is an EP
operator, which according to Remark 12, is equivalent to the 
fact that $a$ is EP.\par 

\indent The third statement can be proved in a similar way
using $R_a$ instead of $L_a$.\par

\indent The last statement is a consequence of Theorem 6.\par
\end{proof}

\indent In the following proposition, the relationship between
EP operators and the adjoint will be studied.\par

\newtheorem*{prop14}{Proposition 14} 
\begin{prop14} Let $X$ be a Banach space, and $T\in L(X)$
such that $T^{\dag}$ exits. Then, necessary and suffcient for
$T$ to be EP is that $T^*\in L(X^*)$ is EP.
\end{prop14}
\begin{proof} Since, according to Proposition 10, $(T^*)^{\dag}=(T^{\dag})^*$,  
it is clear that if $T$ is EP, then $T^*$ is EP.\par

\indent On the other hand, if $T^*$ is EP, then, according
to what has been proved, $(T^*)^*$ also is.
Consider, as in Proposition 10, the isometric isomorphism 
$\Lambda\colon X\to \Lambda (X)\subseteq (X^*)^*$ and $\tilde{T}$, the restriction of
$(T^*)^*$ to $\Lambda(X)$. Then, according to the last paragraph of Remark 12,
$\tilde{T}$ is EP. Now well, since $\tilde{T}^{\dag}$
is the restriction of $((T^{\dag})^*)^*$ to $\Lambda (X)$ and
$$
\Lambda TT^{\dag}= \tilde{T}\tilde{T}^{\dag}\Lambda, \hskip1cm 
\Lambda T^{\dag}T=\tilde{T}^{\dag} \tilde{T}\Lambda,
$$
$TT^{\dag}=T^{\dag}T$.\par
\end{proof}
\indent In [6] a well-known characterization of EP Hilbert
space operators was stated. This result was firstly extended to
$C^*$-algebras in [12], and in second place in [15], where other equivalent
statements were proved and the main concept used was the Drazin inverse. 
The most important results of this
section consist in the extension and reformulation of the aforesaid characterizations
in Banach spaces and 
algebras, using instead of the adjoint, the Moore-Penrose inverse and
the properties of hermitian projectors as developed in 
[3] and [20]. In order to begin with this subject,
the characterization of [6] is considered.\par 
    
\newtheorem*{rem15}{Remark 15} 
\begin{rem15}\rm Let $H$ be a Hilbert space and
let $A\in L(H)$ be an operator with closed range. Recall that if $A^{\dag}$ is the Moore-Penrose of $A$,
then $N(A^{\dag})= N(A^*)$ and $R(A^{\dag})= R(A^*)$. 
Moreover, considering the $C^*$-algebra $L(H)$, according to the proof of Theoreme 10 in [12], 
see Remark 19, there are invertible operators $P_1$ and $P_2$ defined on $H$
such that $A^*= P_1A^{\dag}= A^{\dag}P_2$.
Therefore, the relevant information contained in the characterization of [6]
consists in the fact that $A$ is EP if and only if
$N(A)=N(A^{\dag})$, $R(A)= R(A^{\dag})$, or $A^{\dag}= \tilde{P}A$,
where $\tilde{P}$ is an invertible operator defined on $H$. In the following theorem, in
the context of Banach spaces,
the characterization of [6] will be reformulated using the range and the null space of
the Moore-Penrose inverse instead of the corresponding spaces of the
adjoint. Furthermore, such characterization will be obtained precisely
 thanks to Theorem 2.2 of [20], which states that hermitian idempotents in Banach spaces
have properties similar to the ones of orthogonal projectors in Hilbert spaces.
 \par 
\end{rem15}

\newtheorem*{theo16}{Theorem 16} 
\begin{theo16}Let $X$ be a Banach space, and consider $T\in L(X)$ such that  
$T^{\dag}$ exists. Then, the following statements
are equivalent:\par
(i) \hskip.27cm $T$ is an EP operator,\par
(ii) \hskip.15cm $N(T)=N(T^{\dag})$,\par
(iii) $R(T)=R(T^{\dag})$,\par
(iv) \hskip.04cm there exists an invertible operator $P\in L(X)$ such that $T^{\dag}=PT$.\par
\end{theo16}
\begin{proof} Recall that since $T^{\dag}$ exists, $TT^{\dag}$ (resp. $T^{\dag}T$)
is a hermitian idempotent such that $R(T)=R(TT^{\dag})$ (resp. $R(T^{\dag})=R(T^{\dag}T$).
Futhermore, according to Theorem 2.2 of [20], $R(T)=R(T^{\dag})$ if and only if
$TT^{\dag}=T^{\dag}T$.\par
\indent Similarly, since $T^{\dag}$ exists, according Remark 2 or to Theorem 4.4 of [9], 
$I-TT^{\dag}$ (resp. $I-T^{\dag}T$) is a
hermitian idempotent such that $N(T^{\dag})= R(I-TT^{\dag})$ (resp. $N(T)= R(I-T^{\dag}T)$).
Consequently, according again to Theorem 2.2 of [20], $T$ is an EP operator
if and only if $N(T)=N(T^{\dag})$.\par
\indent Next suppose that there exists an invertible operator $P\in L(X)$ such that
 $T^{\dag}=PT$. It is clear that $N(T)=N(T^{\dag})$. In particular, $T$
is an EP operator.\par
\indent On the other hand, if $T$ is an EP operator, it is not difficult to prove that
$X=N(T)\oplus R(T)$, and $ \tilde{T}= T\mid_{R(T)}\colon R(T)\to R(T)$
is an invertible operator. Define the operator $P'$ in the 
following way:\par
$$
P'\mid_{N(T)}\equiv I_ {N(T)},\hskip1cm P'\mid_{R(T)}= \tilde{T}^2,
$$ 
where $I_ {N(T)}$ denotes the identity map on $N(T)$.\par
\indent Clearly, $P'$ is invertible, and a straightforward calculation
shows that  $T=P'T^{\dag}$. In order to conclude the proof, define $P=(P')^{-1}$.\par
\end{proof}

\newtheorem*{rem17}{Remark 17} 
\begin{rem17}\rm In the same conditions of Theorem 16, note that a
straightforward calculation proves the following identities: $T^{\dag}=PT=TP$.\par
\indent On the other hand, observe that $T$ is an EP operator if and only if
there exists an invertible operator $Q\in L(X)$ such that  $T=QT^{\dag}=T^{\dag}Q$.
In fact, if such identites hold, then $N(T)=N(T^{\dag})$. Consequently, according to
the second statement of Theorem 16, $T$ is an EP operator.
In order to prove the converse implication, note that $T$ is an EP operator
if and only if $T^{\dag}$ is. Then, according to the last statement of Theorem 16
and what has been proved, 
there is an invertible operator $Q$ such that $T=QT^{\dag}=T^{\dag}Q$.\par
\indent Recall that in the Hilbert space context, the operator $A$ has a Moore-Penrose 
inverse if and only if $A^*$, $A^{\dag}$ and $(A^{\dag})^*$ also have one. Since in [6]
and in Theorem 16 the condition of being EP has been characterized using a relationship
between $A$ and $A^*$ and between $A$ and $A^{\dag}$, more equivalent statements
can be deduced considering $A^*$, $A^{\dag}$, and $(A^{\dag})^*$, and their corresponding
adjoints or Moore-Penrose inverses. \par
\indent Finally, note that for any closed range operator $A$ defined on a Hilbert space,
$N(A)=N((A^{\dag})^*)$ and $R(A)=R((A^{\dag})^*)$. Therefore, no characterization of the EP
condition involving $A$ and $(A^{\dag})^*$ can be obtained.\par
\end{rem17}

\indent Next follows the characterization of EP Banach algebra elements.
As in Theorem 16, using the left and right multiplication
representations and ideas of [12] and [15], this result reformulates and extends Theorem 10 of [12]
and Theorem 3.1 of [15] to Banach algebras.\par

\newtheorem*{theo18}{Theorem 18} 
\begin{theo18} Let $A$ be a unital Banach algebra, and consider
$a\in A$ such that $a^{\dag}$ exists. Then, the following statements are equivalent:\par
(i) $a$ is an EP Banach alebra element,\par
(ii) $L_a\in L(A)$ is an EP operator,\par
(iii) $a^{-1}(0)=(a^{\dag})^{-1}(0)$,\par
(iv) $aA=a^{\dag}A$,\par
(v) there exists an invertible operator $P\in L(A)$ such that $L_{a^{\dag}}=$\par
 $PL_a=L_aP$,\par
(vi) there exists an invertible operator $Q\in L(A)$ such that $L_a=$\par
$QL_{a^{\dag}}=L_{a^{\dag}}Q$,\par
(vii) $R_a\in L(A)$ is an EP operator,\par
(viii) $a_{-1}(0)=(a^{\dag})_{-1}(0)$,\par
(ix) $Aa=Aa^{\dag}$,\par
(x) there exists an invertible operator $U\in L(A)$ such that $R_{a^{\dag}}=$\par
$UL_a=R_aU$,\par
(xi) there exists an invertible operator $V\in L(A)$ such that $R_a=$\par
$VR_{a^{\dag}}=R_{a^{\dag}}V$,\par
(xii) $a^2a^{\dag}= a=a^{\dag}a^2$,\par
(xiii) $a\in a^{\dag}A\cap Aa^{\dag}$,\par
(xiv) $a^{\dag}\in aA\cap Aa$,\par
(xv) $a\in a^{\dag}A^{-1}\cap A^{-1}a^{\dag}$,\par
(xvi) $a^{\dag}\in aA^{-1}\cap A^{-1}a$,\par
(xvii) $aA^{-1}=a^{\dag}A^{-1}$,\par
(xviii) $A^{-1}a=A^{-1}a^{\dag}$.\par
\end{theo18}  
\begin{proof} Since $A$ is a unital Banach algebra, according to
Theorem 5 and Remark 12, the first statement is equivalent to the
second and the seventh. \par

\indent In addition, according to Theorem 16, the second, third, fourth, fifth and sixth statements
are equivalent. Furthermore, according to the same theorem, the seventh,
eighth, ninth, tenth and eleventh statements are equivalent.\par

\indent On the other hand, recovering an argument used in Theorem 10 of [12],
it is possible to prove that the first and the twelfth statements are
equivalent.\par 

\indent The equivalence between the twelfth and the  thirteenth statements
can be proved as in Theorem 3.1 of [15]. In addition,
in order to prove that the first and the fourteenth statements are 
equivalent, use what has been proved and the fact that $a$ is EP if and only if $a^{\dag}$ is.\par

\indent Recovering another argument used in Theorem 10 of [12],
it is clear that the first statement implies the fifteenth.
On the other hand, the later statement implies the third and the eighth.
Moreover, in order to prove that the first and the sixteenth 
statements are equivalent, use what has been proved and the fact 
that $a$ is EP if and only if $a^{\dag}$ is.\par

\indent If $a$ is an EP element, then according to the proof of Theorem 10 of [12],
there exist two invertible elements $c$ and $d$, such that $a^{\dag}=ac$ and $a=a^{\dag}d$,
which, according to a straightforward calculation proves that $aA^{-1}=a^{\dag}A^{-1}$.
On the other hand, if the seventeenth statement holds, then it is not difficult to
prove that $a_{-1}(0)=(a^{\dag})_{-1}(0)$. \par

\indent Finally, if $a$ is an EP element, then according again to the proof of Theorem 10 of [12],
the invertible elements of the previous paragraph $c$ and $d$ also satisfy that $a^{\dag}=ca$ and $a=da^{\dag}$,
which, according to a straightforward calculation proves that $A^{-1}a=A^{-1}a^{\dag}$.
On the other hand, if the last statement holds, then it is not difficult to
prove that $a^{-1}(0)=(a^{\dag})^{-1}(0)$. \par

\end{proof}

\newtheorem*{rem19}{Remark 19} 
\begin{rem19}\rm Let $A$ be a $C^*$-algebra and $a\in A$. Recall that according to 
Theorem 10 of [12], if $a^{\dag}$ exists, then
$$
a^*\in A^{-1}a^{\dag} \hskip.5cm \hbox{ and}\hskip.5cm a^*\in a^{\dag} A^{-1}.
$$
\indent According to these facts, it is not difficult to prove that

\begin{align*} 
 & (a^*)^{-1} (0)= (a^{\dag})^{-1} (0),&  &a^*_{-1} (0)= a^{\dag}_{-1} (0), &\\ 
&a^*A=a^{\dag}A, & &Aa^*=Aa^{\dag},&\\
&a^*A^{-1}=a^{\dag}A^{-1}, & &A^{-1}a^*=A^{-1}a^{\dag}.
\end{align*}
\indent Consequently, as in the case of a closed range Hilbert space
operator, the relevant equivalences in Theorem 10 of [12]
and Theorem 3.1 of [15] are the following identities:
\begin{align*} 
 & a^{-1} (0)= (a^{\dag})^{-1} (0),&  &a_{-1} (0)= a^{\dag}_{-1} (0), &\\ 
&aA=a^{\dag}A, & &Aa=Aa^{\dag},&\\
&aA^{-1}=a^{\dag}A^{-1}, & &A^{-1}a=A^{-1}a^{\dag}.
\end{align*}
\indent On the other hand, in Theorem 18,
Theorem 10 of [12] and
Theorem 3.1 of [15] it has been characterized the
condition of being EP using a relationship between 
$a$ and $a^{\dag}$ or between $a$ and $a^*$.
Since $a\in A$ is an EP element
if and only if $a^*$, $a^{\dag}$, and $(a^{\dag})^*$ also are,
it is possible to obtain more
equivalent statements for an element to be EP,
applying the mentioned results to $a^*$, $a^{\dag}$, and $(a^{\dag})^*$,
and to their corresponding adjoints and Moore-Penrose inverse.\par

\indent Observe that, as in the case of a closed 
range Hilbert space operator, according to Proposition 2.4 of [4],
it is not difficult to prove that 

\begin{align*} 
 & a^{-1} (0)= ((a^{\dag})^*)^{-1} (0),&  &a_{-1} (0)= (a^{\dag})^*_{-1} (0), &\\ 
&aA=(a^{\dag})^*A, & &Aa=A(a^{\dag})^*.
\end{align*}

\indent Therefore, in order to characterize EP elements,
no characterization involving $a$ and $(a^{\dag})^*$ can be
obtained.\par

\indent Finally, if $X$ is a Banach space and $T\in L(X)$
is an operator with a Moore-Penrose inverse, Theorem 18
provides a set of characterizations for $T$ to be EP,
considering it as an element of the Banach algebra $A=L(X)$.\par 
\end{rem19}

\noindent \bf{5. The Product of Two EP Elements }\rm \vskip.3cm

\indent In the present section it will be studied
when the product of two EP Banach space operators or
Banach algebra elements is again EP. In first place,
a remark is considered.\par

\newtheorem*{rem20}{Remark 20} 
\begin{rem20}\rm The problem of characterizing 
when the product of two EP matrices is again EP 
was first posed in [1], and solved in Theorem 1 of [13], where 
it was formulated using the row space of a matrix. In 
[14] a simple proof of the latter Theorem was given
using an operator theoretical language. \par
\indent For closed range Hilbert space operators, in [8] it was proved that
the first and the fourth statements of Theorem 2 of [14] are equivalent, while the Example in section 3,
page 115, of [16]
shows that the third and 
the fourth statements of Theorem 2 of [14] are not any more equivalent. On the other hand, in Theorem 4.3 and Corollary
4.4 of [15] the problem under consideration was studied in the context of
$C^*$-algebras and closed range Hilbert space operators. In the following 
theorems, using hemitian idempotents and reformulating an idea of
Theorem 4.3 of [15], it will be characterized when the product of two
Banach algebra elements or Banach space operators is again EP.\par
\indent Two more observations. Recall that if $X$ is a Banach space, then an operator
$T\in L(X)$ is called \it upper semi-Fredholm \rm (resp. \it lower semi-Fredholm), \rm
if $R(T)$ is closed and $N(T)$ is finite dimensional (resp. $R(T)$
has finite codimension). If $T\in L(X)$ is EP,
then necessary and sufficient for $T$ to be a Fredholm operator is that
$T$ is upper semi-Fredholm (resp. lower semi-Fredholm).
Moreover, in this case $ind (T)=0$, where $ind$ denotes
the index of $T$. On the other hand, if $T\in L(X)$ is an operator
such that $T^{\dag}$ exists, then $I-T^{\dag}T$ is the 
hermitian idempotent onto $N(T)$, see Theorem 2 of [14]. \par  
\end{rem20}

\newtheorem*{theo21}{Theorem 21} 
\begin{theo21}Let $X$ be a Banach space, and let $S$ and
$T$ be two EP Banach space operators defined on $X$
such that the Moore-Penrose inverse of $ST$ exists.\par
(i) If $N(ST)=N(S) +N(T)$ and $R(ST)=R(S)\cap R(T)$, then
$ST$ is  EP.\par
(ii) If $ST$ is an EP Banach space operator, then $R(ST)\subseteq R(T)$
and $N(S)\subseteq N(ST)$.\par
(iii) Necessary and sufficient for $R(ST)\subseteq R(T)$ (resp. $N(S)\subseteq N(ST)$) is that
$(I-TT^{\dag})ST=0$ (resp. $ST(I-S^{\dag}S)=0$).\par
(iv) If $(I-T^{\dag}T)ST=0$, $ST(I-S^{\dag}S)=0$,
and $S$ and $T$ are Fredholm operators, then $N(ST)=N(S) +N(T)$ and $R(ST)=R(S)\cap R(T)$.
In particular, $ST$ is EP.   
\end{theo21}

\begin{proof} Suppose that $N(ST)=N(S) +N(T)$ and $R(ST)=R(S)\cap R(T)$.
Define $M=N(S)$ and $N=N(T)$. Then, $M$ and $N$ belong to
$\mathcal{M}(X)$, see Remark 2. In fact, $P_M=I-S^{\dag}S$ and 
$P_N=I-T^{\dag}T$. In addition, $M'=N(P_M)=R(S)$ and $N'=N(P_N)=R(T)$.\par
\indent Now well, since $ST(ST)^{\dag}$ is a hermitian idempotent such that
$$
R(ST(ST)^{\dag})=R(ST)=R(S)\cap R(T)=M'\cap N',
$$
$M'\cap N'\in \mathcal{M}(X)$ and $P_{M'\cap N'}=ST(ST)^{\dag}$. \par
\indent On the other hand, since $(ST)^{\dag}ST$ is an hermitian idempotent
such that 
$$
R(I-(ST)^{\dag}ST)=N(ST)=N(S)+N(T)=M+N,
$$
$M+N$ belongs to $\mathcal{M}(X)$ and $P_{M+N}=I-(ST)^{\dag}ST$.\par
\indent Therefore, according to Theorem 2.28 of [3],
$$
P_{\overline{M+N}} =I-P_{M'\cap N'}.
$$
\indent However, $\overline{M+N}=N(ST)$ and $M'\cap N'=R(ST)$.
Consequently,
$$
I-(ST)^{\dag}ST=I-ST(ST)^{\dag},
$$
equivalently, $(ST)^{\dag}ST=ST(ST)^{\dag}$.\par
 
\indent Next suppose that $ST$ is an EP Banach space operator.
Define $M_1= R(ST)$ and $N_1=R(T)$. Then $M_1$ and $N_1$
belong to $\mathcal{M}(X)$, $P_{M_1}=ST(ST)^{\dag}$, $P_{N_1}=TT^{\dag}$,
$M_1'=N((ST)^{\dag})$ and $N_1'=N(T^{\dag})$.\par
\indent Now well, according to the observation that precedes
Theorem 2.28 of [3], equivalent for $R(ST)$ to be contained 
in $R(T)$ is that $N(T^{\dag})$ is containd in $N((ST)^{\dag})$.
However, since $T$ and $ST$ are EP Banach space operators,
according to Theorem 16, $N(T^{\dag})=N(T)\subseteq N(ST)=
N((ST)^{\dag})$.\par
\indent On the other hand, consider $M_2=R((ST)^{\dag})$ and $N_2=R(S)$. 
As before, $M_2$ and $N_2$ belong to $\mathcal{M}(X)$, 
$P_{M_2}=(ST)^{\dag}ST$, $P_{N_2}=SS^{\dag}$, $M_2'=N(ST)$,
and $N_2'=N(S^{\dag})$.\par
\indent Now well, according again to the observation that precedes Theorem
2.28 of [3], necessary and sufficient for $N(S)$ to be contained in 
$N(ST)$ is that $R((ST)^{\dag})$ is contained in $R(S^{\dag})$.
However, due to the fact that $S$ and $ST$ are EP Banach space operators,
according to Theorem 16, $R((ST)^{\dag})=R(ST)\subseteq R(S)=R(S^{\dag})$.\par
\indent As regard the third statement, in order to prove the equivalence it is
enough to observe that $R(T)=R(TT^{\dag})=N(I-TT^{\dag})$ and $N(S)=N(S^{\dag}S)=
R(I-S^{\dag}S)$.\par

\indent In order to prove the last statement, suppose that 
$$
ST(I-S^{\dag}S)=0 \hbox{ and } (I-TT^{\dag})ST=0.
$$
\indent Next observe that, since $S$ and $T$ are EP Banach space
operators, 
$$
(I-S^{\dag}S)ST=0 \hbox{ and } ST(I-TT^{\dag})=0.
$$
Therefore,
$$
STS^{\dag}S=S^{\dag}SST=ST \hbox{ and } STTT^{\dag}=TT^{\dag}ST=ST.
$$
\indent Conseqently, $R(ST)$ and $N(ST)$ are closed invariant subspaces for
$S^{\dag}S$, $TT^{\dag}$, $I-S^{\dag}S$, and $I-TT^{\dag}$. In particular,
according to Proposition 4.12 of [9], $TT^{\dag}$ and $I-TT^{\dag}$
are hermitian idempotents when restricted to $R(ST)$ and $N(ST)$,
and due to the fact that $(I-TT^{\dag})(X)=(I-TT^{\dag})(N(ST))=N(T)$, it is clear that
$$
R(ST)=TT^{\dag} (R(ST)) \hbox{ and } N(ST)= TT^{\dag} (N(ST))\oplus N(T).
$$ 
\indent Define the linear and bounded operator $U\in L(N(ST))$ in the following way:
$U\mid N(T)=I_{N(T)}$, the identity map of $N(T)$, and 
$U=T\mid_{TT^{\dag}(N(ST))}^{T(N(ST))}\colon TT^{\dag}(N(ST))\to T(N(ST))$. Note that $T(N(ST))\subseteq R(T)\cap N(S)\subseteq N(S)
\subseteq N(ST)$. Moreover, a straightforward calculation proves that $T(N(ST))=R(T)\cap N(S)$. 
What is more, since $TT^{\dag}(N(ST))\subseteq R(T)$ and $N(T)\cap R(T)=0$,
$U\colon N(ST)\to R(T)\cap N(S)\oplus N(T)$ is a Banach space isomorphism.\par
\indent Now well, if $N(S)$ and $N(T)$ are finite dimensional, then $R(U)=N(ST)= R(T)\cap N(S)\oplus N(T)$. 
Therefore, $N(ST)=N(S)+N(T)$.\par

\indent In order to prove that $R(ST)=R(S)\cap R(T)$, a duality argument will be
considered. \par
\indent Observe that, according to the hypothesis of the last statement of the Theorem,
$S^*$ and $T^*$ are Fredholm operators defined on $X^*$, the dual space
of $X$.  Moreover, according to Proposition 14 and Theorem 10,
$S^*$ and $T^*$ are EP, and $T^*S^*$ has a Moore-Penrose inverse. In addition, it is clear that 
$$
(I-S^*(S^*)^{\dag})T^*S^*=0 \hbox{ and } T^*S^*(I-(T^*)^{\dag}T^*)=0.
$$
\indent Therefore, according to what has been proved, $N(T^*S^*)=N(T^*)+N(S^*)$.
Moreover, since $ST$ is a Fredholm operator, its range is a closed subspace of $X$ and
\begin{align*} 
R(ST)&=\hbox{ }^{\perp}N((ST)^*)=\hbox{ }^{\perp}(N(S^*)+N(T^*))=\hbox{ }^{\perp}N(S^*)\hbox{ }\cap
\hbox{ }^{\perp}N(T^*)\\
&=R(S)\cap R(T),
\end{align*} 
where if $V$ is a subspace of $X^*$, then $\hbox{ }^{\perp}V=\{ x\in X\colon f(x)=0,\hbox{ } \forall f\in V \}$.
\end{proof}

\indent In the following Theorem it will be studied when the product of two Banach algebra
elements is EP.\par

\newtheorem*{theo22}{Theorem 22} 
\begin{theo22} Let $A$ be a unital Banach algebra. Consider $a$ and $b$,
two EP elements of $A$ such that $ab$ has a Moore-Penrose inverse.\par
(i) If $(ab)^{-1}(0)=a^{-1}(0)  +b^{-1}(0)$ and $abA=aA\cap bA$, then $ab$ is EP.\par
(ii) If $ab$ is an EP Banach algebra element, then $abA\subseteq bA$ and 
$a^{-1}(0)\subseteq (ab)^{-1}(0)$.\par
(iii) Necessary and sufficient for $abA\subseteq bA$ (resp. $a^{-1}(0)\subseteq  (ab)^{-1}(0)$) 
is that $(e-bb^{\dag})ab=0$ (resp. $ab(e-a^{\dag}a)=0$).\par
(iv)  If $(e-bb^{\dag})ab=0$, $ab(e-a^{\dag}a)=0$, and  $a^{-1}(0)$ and 
$b^{-1}(0)$ are finite dimensional, then $(ab)^{-1}(0)=a^{-1}(0)  +b^{-1}(0)$ and $abA=aA\cap bA$.
In particular, $ab$ is EP.\par
\indent Similarly,\par
(v) If $(ab)_{-1}(0)=a_{-1}(0)  +b_{-1}(0)$ and $Aab=Aa\cap Ab$, then $ab$ is EP.\par
(vi) If $ab$ is an EP Banach algebra element, then $Aab\subseteq Aa$ and 
$b_{-1}(0)\subseteq (ab)_{-1}(0)$.\par
(vii) Necessary and sufficient for $Aab\subseteq Aa$ (resp. $b_{-1}(0)\subseteq (ab)_{-1}(0)$) 
is that $ab(e-aa^{\dag})=0$ (resp. $(e-b^{\dag}b)ab=0$).\par
(viii)  If $ab(e-aa^{\dag})=0$, $(e-b^{\dag}b)ab=0$, and $a_{-1}(0)$ and 
$b_{-1}(0)$ are finite dimensional, then $(ab)_{-1}(0)=a_{-1}(0)  +b_{-1}(0)$ and $Aab=Aa\cap Ab$.
In particular, $ab$ is EP.\par

\begin{proof} Consider $L_a$ (resp. $R_a$), $L_b$ (resp. $R_b$), and $L_{ab}$
(resp. $R_{ab}$), and then apply Theorem 21, and use Theorem 5 and Remark 12.\par
\end{proof}
\end{theo22} 

\noindent \bf{Acknowledgements.} \rm The author wish to express his indebtedness to 
Professor Valdimir Rakocevic, for this researcher has 
turned the attention of the author to the Moore-Penrose
inverse in Banach algebras, and to the referee, for several suggestions that
have improved the final version of the present work.\par

\vskip.3truecm
\noindent Enrico Boasso\par
\noindent E-mail address: enrico\_odisseo@yahoo.it

\end{document}